\documentclass[11pt]{amsart}

\begin{document}

\newtheorem{lemma}{Lemma}[section] 
\newtheorem{theorem}{Theorem}[section]
\newtheorem{definition}{Definition}[section]
\newtheorem{proposition}{Proposition}[section]
\newtheorem{corollary}{Corollary}[section]
\newtheorem{conjecture}{Conjecture}[section]

\title{Growth estimates for discrete quantum groups} 

\author{Teodor Banica}
\address{T.B.: Department of Mathematics, Toulouse 3 University, 118 route de Narbonne, 31062 Toulouse, France}
\email{teodor.banica@math.ups-tlse.fr}

\author{Roland Vergnioux} 
\address{R.V.: Department of Mathematics, Caen University, BP 5186, 14032 Caen Cedex, France}
\email{roland.vergnioux@math.unicaen.fr} 


\begin{abstract}
We discuss the notion of growth for discrete quantum groups, with a number of general considerations, and with some explicit computations. Of particular interest is the quantum analogue of Gromov's estimate regarding  polynomial growth: we formulate the precise question, and we verify it for the duals of classical Lie groups.
\end{abstract}

\maketitle

\section*{Introduction}

The present paper is motivated by some basic questions regarding the compact quantum groups.  These are abstract objects, corresponding to the Hopf algebras in the sense of Woronowicz. We consider compact quantum groups which are coinvolutive, or of Kac type, in the sense that they satisfy one of the following equivalent conditions:
\begin{enumerate}
\item The square of the antipode is the identity. 

\item The Haar functional is a trace.
\end{enumerate}

Such a Hopf algebra $A$ produces a compact quantum group $G$ and a discrete quantum group $\Gamma$, according to the following heuristic formula:
$$A={\rm C}(G)={\rm C}^*(\Gamma)$$
In other words, studying compact quantum groups is the same as studying discrete quantum groups. This point of view is essential for considerations in this paper, where ``compact'' and ``discrete'' methods are used at the same time.

The motivations for compact quantum groups come in part from a number of algebraic and analytic problems, explained below, and on the other hand, from an intimate link with certain combinatorial aspects of von Neumann algebras, free probability, statistical mechanical models, or random matrices, discussed in \cite{b3}, \cite{bbc}, \cite{bc}.

So, let $G$ be a compact quantum group. The finite dimensional unitary representations of $G$ form a fusion semiring $R^+$, coming with a dimension function $\dim:R^+\to {\mathbb N}$. There are some natural problems regarding the invariant $(R^+,\dim)$, discussed some time ago in the survey paper \cite{b1}. Here is an updated list:
\begin{enumerate}
\item Deformation conjecture. This states that for any compact quantum group, there are $1\leq n<\infty$ compact quantum groups of Kac type having the same $(R^+,\dim)$ invariant. The ``anti-deformation'' inequality $n\geq 1$ was checked in many situations, see \cite{b1}. As for the finiteness inequality $n<\infty$, this holds in the finite quantum group case, by a result of Stefan \cite{st}.

\item Operator algebras. Free quantum groups provide new examples of simple ${\rm C}^*$-algebras, K-amenable quantum groups, and prime ${\rm II}_1$ factors, see \cite{b1}, \cite{vv}, \cite{v1}. The proofs so far are quite technical, but seem to ultimately rely on special properties of the corresponding $(R^+,\dim)$ invariants.

\item Polynomial growth. A fundamental result of Gromov states that a finitely generated group has polynomial growth if and only if it is virtually nilpotent \cite{gr}. For quantum groups such subtleties are out of sight for the moment: in fact, polynomial growth itself is quite tricky to define, in terms of $(R^+,\dim)$, and this was worked out only recently, in \cite{v3}.

\item Tannakian philosophy. A fundamental result of Deligne \cite{de} and Doplicher-Roberts \cite{dr} states that a symmetric tensor category with duals comes from a compact group if and only if all internal dimensions are integers. In the quantum group case there are many difficult problems around, for instance characterization of subfactors of integer index, or of values of $(R^+,\dim)$. See \cite{b1}.
\end{enumerate}

In fact, very little is known about $(R^+,\dim)$. However, there are some general directions to explore, coming from both compact and discrete groups:

\begin{enumerate}
\item Amenability. A lower bound for dimensions for amenable groups comes from work of Longo and Roberts in \cite{lr}, and from the Kesten amenability type criterion in \cite{b2}. This provides the only known general restriction on $(R^+,\dim)$.

\item Weyl formula. The problem with compact quantum groups is that no analogue of the Lie algebra is available. However, it seems reasonable to think that the Weyl formula should have extensions to certain quantum groups.  But, as far as we know, there is no precise result of that type.

\item Poincar\'e duality. For $R$-matrix quantizations of $SU(n)$ such a duality is obtained by Gurevich in \cite{gu}. This seems to impose restrictions on $(R^+,\dim)$, but once again, to our knowledge, there is no further result of that type.

\item Growth. A number of key estimates are known to relate growth and random walks on discrete groups, see the book \cite{vsc}.  These provide subtle analytic restrictions on the invariant $(R^+,\dim)$ for duals of discrete groups.
\end{enumerate}

Summarizing, we have a number of difficult problems regarding compact quantum groups, and a ray of light coming from the concept of growth of discrete groups.

The purpose of this paper is to make some advances on the subject. We first review the notion of growth for discrete quantum groups, recently introduced in \cite{v3}. Then we perform a few explicit computations: for various product operations, for free quantum groups, and for duals of Lie groups. The free quantum groups are studied by using linearization formulae for free product operations, and have exponential growth. The duals of Lie groups are known from \cite{v3} to have polynomial growth, and we compute here the growth exponent, by improving the way the Weyl formula is applied.

The second part of the paper discusses the notion of random walk.  The relevant quantity here is the probability of returning to the origin, after $k$ steps on the Cayley graph. This can be in turn expressed in terms of representation theory invariants, such as multiplicities, or integrals of characters. For free quantum groups this kind of quantity is well-known, and we write down the explicit first order estimates. For duals of Lie groups the computation is quite technical, and once again by using the Weyl formula, along with some analytic techniques, we compute the exponent.

The computations of growth and random walks that we have turn out to be related by estimates which are similar to those in \cite{vsc}.  Of particular interest here is the quantum group formulation of a key estimate of Gromov, regarding polynomial growth. This quantum group statement is a conjecture, and the evidence comes from Gromov's result, and from the Lie group computations in this paper. In other words, the conjecture is basically verified in both the classical and the dual classical case. Previous experience with discrete quantum groups tells us that in this situation, the conjecture should indeed be true. However, both verifications are quite technical, and we don't see yet the connection between them. In other words, we are still far away from a possible proof, and the whole paper should be regarded as an introduction to the subject.

Finally, let us mention that the Lie group considerations in this paper are closely related to the machinery developed in \cite{bia}, \cite{bbo}, \cite{co}. Also, the idea of investigating nilpotent quantum groups appears at level of tensor categories in \cite{gn}.

The paper is organized as follows: 1 is a preliminary section, in 2-4 we compute growth invariants, and in 5 we discuss random walks and we state the conjecture.

\section{Discrete quantum groups}

Let $A$ be a unital ${\rm C}^*$-algebra, given with a morphism $\Delta :A\to A\otimes A$ satisfying the coassociativity condition $(\Delta\otimes {\rm id})\Delta =({\rm id}\otimes\Delta )\Delta$, and with a morphism $\varepsilon: A\to{\mathbb C}$ and an antimorphism $S:A\to A^{op}$ satisfying the usual axioms for a counit and for an antipode, along with the coinvolutivity condition $S^2={\rm id}$. For the precise formulation of the axioms see Woronowicz \cite{wo}.

For simplicity of presentation, we use the following terminology:

\begin{definition}
In this paper we call Hopf algebra a unital ${\rm C}^*$-algebra $A$, given with a comultiplication $\Delta$, a counit $\varepsilon$, and an antipode $S$ satisfying $S^2={\rm id}$.
\end{definition}

The first example is related to a discrete group $\Gamma$. We can consider the group algebra $A={\rm C}^*(\Gamma)$, with operations given by
\begin{eqnarray*}
\Delta(g)&=&g\otimes g\\
\varepsilon(g)&=&1\\
S(g)&=&g^{-1}
\end{eqnarray*}
where we use the canonical embedding $\Gamma\subset A$.

The second example is related to a compact group $G$. We have here the algebra of continuous functions $A={\rm C}(G)$, with operations given by
\begin{eqnarray*}
\Delta(\varphi)&=&(g,h)\to \varphi(gh)\\
\varepsilon(\varphi)&=&\varphi(1)\\
S(\varphi)&=&g\to \varphi(g^{-1})
\end{eqnarray*}
where we use the canonical identification $A\otimes A={\rm C} (G\times G)$.

In general, associated to a Hopf algebra $A$ are a discrete quantum group $\Gamma$ and a compact quantum group $G$, according to the following heuristic formulae:
$$A={\rm C}^*(\Gamma)={\rm C}(G).$$
For understanding what $\Gamma$ is, it is best to use representations of $G$. These correspond to corepresentations of $A$, which can be defined in the following way.

\begin{definition}
A corepresentation of $A$ is a unitary matrix $u\in M_n(A)$ satisfying
\begin{eqnarray*}
\Delta(u_{ij})&=&\sum u_{ik}\otimes u_{kj}\\
\varepsilon(u_{ij})&=&\delta_{ij}\\
S(u_{ij})&=&u_{ji}^*
\end{eqnarray*}
for any $i,j$. We say that $u$ is faithful if its coefficients generate $A$ as a ${\rm C}^*$-algebra. A finitely generated Hopf algebra $(A,u)$ is a Hopf algebra $A$ together with a faithful corepresentation $u$ such that $u=\bar{u}$.
\end{definition}

There are notions of direct sum, tensor product, complex conjugation, equivalence etc., analogous to those for representations of compact groups. We use the notations $u+v,u\otimes v,\bar{u}$ for sums, tensor products, and complex conjugation. The condition $u=\bar{u}$ can be achieved by replacing $u$ with $u+\bar{u}$ if necessary.

For the algebra $A={\rm C}(G)$, the corepresentations of $A$ are in one-to-one correspondence with the unitary representations of $G$ via the canonical identifications $M_n(A)={\rm C}(G, M_n({\mathbb C}))$. Moreover the operations $+$, $\otimes$, $\bar~$ correspond to the direct sum, tensor product and contragredience operations, and faithfulness is the usual notion. 

In particular finitely generated Hopf algebras $({\rm C}(G),u)$ correspond to faithfully represented compact Lie groups. In this paper we will mainly deal with connected simply connected compact real Lie groups, and we will endow such a Lie group $G$ with the direct sum $u$ of its fundamental representations -- by this we mean the irreducible representations associated to the fundamental weights of a root system of $G$, which are well-defined up to equivalence.

For the algebra $A={\rm C}^*(\Gamma)$, the corepresentations of $A$ decompose as finite sums of one-dimensional corepresentations, which in turn are given by elements of $\Gamma$. More precisely, a discrete group $\Gamma$ can be recovered from its Hopf algebra $A={\rm C}^*(\Gamma)$ as stated in the following result:
\begin{enumerate}
\item $\Gamma$ is the set of (classes of) irreducible corepresentations of $A$.

\item The product of $\Gamma$ is the tensor product of corepresentations $\otimes$.

\item The unit of $\Gamma$ is the trivial corepresentation $1$.

\item The inverse operation of $\Gamma$ is the complex conjugation.
\end{enumerate}
In particular any corepresentation $u$ must be a direct sum of elements of $\Gamma$, in the corepresentation theory sense. In other words we have
$$u=\begin{pmatrix}g_1&&0\cr&\ddots&\cr0&&g_n\end{pmatrix}$$
for some elements $g_1,\ldots ,g_n\in\Gamma$. Now the fact that the coefficients of $u$ generate $A$ translates into the fact that the set $U=\{g_1,\ldots,g_n\}$ generates $\Gamma$, so that finitely generated Hopf algebras $({\rm C}^*(\Gamma),u)$ correspond to finitely generated groups $(\Gamma,U)$, where $U$ is stable by inversion.

These remarks tell us what a discrete quantum group is: we have a set $\Gamma$, all whose elements are represented by matrix blocks $M_n({\mathbb C})$. The analogues of the product, of the unit element, and of the inverse operation are
the fusion rules
$$u\otimes v=w_1+w_2+\ldots +w_s,$$
the trivial corepresentation $1$, and the complex conjugation operation $u\to\bar{u}$. Here $\Gamma$ is of course the set of irreducible corepresentations of $A$, and the matrix block $M_n({\mathbb C})$ representing a corepresentation $v$ is given by $n=\dim(v)$.

We have now all ingredients for a key definition, first used in \cite{v3} and which generalises the usual notions for discrete groups:

\begin{definition}
Let $(A,u)$ be a finitely generated Hopf algebra.
\begin{enumerate}
\item $\Gamma$ is the set of irreducible corepresentations of $A$ (up to equivalence).

\item The length of $v\in\Gamma$ is the number $l(v)=\min\{ k\in{\mathbb N}\mid v\subset u^{\otimes k}\}$.

\item The volume of balls and spheres are given by
$$b_k=\sum_{l(v)\leq k} \dim(v)^2 \quad\text{and}\quad s_k=\sum_{l(v) = k} \dim(v)^2.$$

\item The $B$ and $S$ series are $B=\sum_{k=0}^\infty b_kz^k$ and $S=\sum_{k=0}^\infty s_kz^k$.
\end{enumerate}
\end{definition}

The fact that lengths are finite follows from Woronowicz's analogue of Peter-Weyl theory \cite{wo}. Actually, one can prove that the formula
$$d(v,w)=\min\{ k\in{\mathbb N}\mid 1\subset \bar{v}\otimes w\otimes u^{\otimes k}\}$$
defines a distance of $\Gamma$, and that we have $l(v)=d(v,1)$. Moreover, the quasi-isometry class of $\Gamma$ doesn't depend on the choice of $u$. See \cite{b1}, \cite{v2}.

In the above definition $z$ is a formal variable, or a complex number where the series converge. We have $s_k=b_k-b_{k-1}$ and $S=B(1-z)$, so the sequences $(b_k),(s_k)$ and the series $B,S$ encode the same information. In practice, $(b_k)$ is most appropriate for asymptotic estimates, whereas $S$ is most convenient for exact formulae.

\begin{definition}
We use the following notations, for sequences of positive numbers:
\begin{enumerate}
\item $b_k\simeq c_k$ means that $b_k/c_k$ converges to a positive number.

\item $b_k\approx c_k$ means $\alpha<b_k/c_k<\beta$, for some constants $0<\alpha,\beta<\infty$.
\end{enumerate}
\end{definition}

\section{Duals of compact groups}

A finitely generated Hopf algebra has polynomial growth when the sequence $b_k$ has polynomial growth, meaning that $b_k$ is dominated by a polynomial in $k$.

The basic example is the algebra $A={\rm C}^*(\Gamma)$, with $\Gamma$ discrete group having polynomial growth. If $b_k \approx k^d$, we call $d$ the exponent of the polynomial growth. This exponent is independent of the chosen faithful corepresentation.

The notion of polynomial growth was introduced in \cite{v3}. In addition to the general results found there, we have the following standard result:

\begin{proposition}
If $A$ has polynomial growth, then $A$ is amenable in the discrete quantum group sense.
\end{proposition}

\begin{proof}
Denoting by $B_k$ the sphere of radius $k$ and center $1$ in $\Gamma$, the tensor powers of $u$ decompose as
$$u^{\otimes k}=\sum_{r\in B_k} m_k(r)\, r$$
where $m_k(r)\in{\mathbb N}$ are certain multiplicities. Due to our assumption $u=\bar{u}$ we have
$$u^{\otimes 2k}=u^{\otimes k}\otimes\bar{u}^{\otimes k}=\sum_{r,p\in B_k} m_k(r)m_k(p)\, r\otimes\bar{p}$$
which by standard results in \cite{wo} gives $m_{2k}(1)=\Sigma_{r\in B_k}m_k(r)^2$. Now the inequality
\begin{eqnarray*}
m_{2k}(1)\times b_k 
&=&\sum_{r\in B_k}m_k(r)^2\times\sum_{r\in B_k}\dim(r)^2\\
&\geq&\Big(\sum_{r\in B_k} m_k(r)\,\dim(r)\Big)^2 = \dim(u)^{2k}
\end{eqnarray*}
shows that if $b_k$ has polynomial growth, then the following happens.
$$\limsup_{k\to\infty}\, m_{2k}(1)^{1/2k}\geq \dim(u)$$
The Kesten type criterion in \cite{b2} applies and gives the result.
\end{proof}

The fact that the duals of compact connected real Lie groups have polynomial growth is known from \cite{v3}. We proceed now to a more precise estimate yielding the growth exponent:

\begin{theorem}
For the discrete quantum group associated to the algebra ${\rm C}(G)$, with $G$ connected simply connected compact real Lie group, we have $b_k\approx k^d$, where $d$ is the real dimension of $G$.
\end{theorem}

\begin{proof}
  Following \cite{bo9}, we fix a maximal torus $T \subset G$, we put $X = \hat
  T$ and $V = X \otimes_{\mathbb Z} {\mathbb R}$. Let $R \subset X$ be the root
  system of $(G,T)$ and endow $V$ with a scalar product invariant under the Weyl
  group of $R$.  Fix a subset of positive roots $R_+ \subset R$ and denote by
  $B\subset R$ the associated basis of $R$. Let $(\varpi_\beta)_{\beta\in B}$ be the
  corresponding basis of fundamental weights, so that we have
  $(\beta'|\varpi_\beta) = 0$ for $\beta \neq \beta'$, and $(\beta | \varpi_\beta) >
  0$. Taking highest weights induces an identification between the set $\Gamma$
  of classes of irreducible representations and the set of dominant weights
  $X_{++}$, which equals in our simply connected case:
  \begin{displaymath}
    X_{++} = \left\{\sum \lambda_\beta \varpi_\beta\, \big\vert\,\lambda_\beta
      \in {\mathbb N}\right\}.
  \end{displaymath}

  Recall that in the case of connected simply connected compact real Lie groups
  we use the word length $l$ associated to the faithful representation $u =
  \sum \varpi_\beta$. Since $\bigotimes \varpi_\beta^{\otimes\lambda_\beta}$
  contains $\sum\lambda_\beta \varpi_\beta$ and other irreducible representations
  with smaller coordinates, we have
  \begin{displaymath}
    l\left(\sum \lambda_\beta \varpi_\beta\right) = \sum \lambda_\beta
  \end{displaymath}
  so $l$ coincides with the restriction to $X_{++}$ of the $\ell^1$-norm
  associated to the basis $(\varpi_\beta)$. Besides, the dimension of $\lambda =
  \sum \lambda_\beta \varpi_\beta$ is given by Weyl's formula
  \begin{displaymath}
    {\dim} (\lambda) = \prod_{\alpha\in R_+} \left( 1 +
      \frac{(\lambda | \alpha)}{(\rho | \alpha)} \right)
  \end{displaymath}
  where $\rho$ is the half-sum of the fundamental weights $\varpi_\beta$.

  We start by estimating dimensions. Since all elements of $R_+$ are sums of
  elements of $B$, one can replace $B$ with $R_+$ in the definition of the
  following closed cones:
  \begin{eqnarray*}
    C&=&\left\{ v\in V ~|~ (v|\alpha)
      \geq 0,\,\forall\, \alpha\in B\right\}\cr
    C'& =& \left\{ v\in V ~|~ (v|\alpha)
      \geq \varepsilon \|v\|,\,\forall\, \alpha\in B\right\}.
  \end{eqnarray*}
  Here $\varepsilon$ is a fixed positive number such that $C'$ is non-empty. By
  definition we have $X_{++} = X \cap C$, and we put $X'_{++} = X \cap C'$.  Let
  $s$ be the number of elements of $R$ and $r={\dim} (T)$, related by the
  formula $s = d - r$. By Weyl's formula and the equivalence of norms on $V$, we
  get positive constants $K_1$, $K_2$ such that
  \begin{eqnarray*}
    {\dim} (\lambda_1)&\leq&K_1 l(\lambda_1)^{s/2}\cr {\dim} (\lambda_2
    )&\geq& K_2 l(\lambda_2)^{s/2}
  \end{eqnarray*}
  for all $\lambda_1\in X_{++}$ and $\lambda_2\in X'_{++}$.

  Let us now estimate the cardinality of the sphere $S_k = \{\lambda\in
  X_{++}\simeq\Gamma ~|~ l(\lambda) = k\}$. By definition, $S_k$ is indexed by
  the following set.
  \begin{displaymath}
    \left\{(\lambda_i)_i \in {\mathbb N}^r\,\big\vert\, \sum\lambda_i = k\right\}
  \end{displaymath}
  It follows that there exist constants $L_1,L_2>0$ such that:
  \begin{displaymath}
    L_1 k^{r-1} \leq \# S_k \leq L_2 k^{r-1}.
  \end{displaymath}
  We need a similar estimate for the cardinality of $S'_k = S_k \cap C'$.  We
  denote by $P_k,P'_k$ the intersection of $C,C'$ with the following hyperplane.
  \begin{displaymath}
    H_k = \left\{(\lambda_i)_i\in {\mathbb R}^r\,\big\vert\,\sum \lambda_i = k \right\}
  \end{displaymath}
  Since $P_k$, $k\in\mathbb{N}$ (resp. $P'_k$, $k\in\mathbb{N}$) are homotetic
  convex polygonal domains, the ratio of volumes
  $\mathrm{vol}(P'_k)/\mathrm{vol}(P_k)$ is constant. By intersecting these
  domains with the lattices $X\cap H_k$ we obtain $S_k$ (resp. $S'_k$), and
  since lattices have fixed covolume we get constants $M_1,M_2>0$ such that:
  \begin{displaymath}
    \# S'_k \geq M_1 (\# S_k) - M_2.
  \end{displaymath}

  We plug now these results into the obvious inequality
  \begin{displaymath}
    (\# S'_k)~ \min_{\lambda\in S'_k} {\dim} (\lambda)^2 \leq
    s_k \leq (\# S_k)~ \max_{\lambda\in S_k} {\dim} (\lambda)^2
  \end{displaymath}
  and obtain the following estimates.
  \begin{displaymath}
    (M_1L_1 k^{r-1} - M_2) K_2^2 k^{s} \leq s_k
    \leq L_2 k^{r-1} K_1^2 k^{s}
  \end{displaymath}
  Now from $s + r =d$ we get $s_k \approx k^{d-1}$, which gives the result.
\end{proof}

We don't know other examples of Hopf algebras having polynomial growth. The problem here is that most work in the area has gone into study of universal examples, which are in general not amenable, hence not of polynomial growth.

A natural idea would be to look at quantum permutation groups.  There are many such quantum groups waiting to be investigated, for instance those which are close to a certain asymptotic area, where non-commutativity is expected to disappear. See \cite{bbc}.

\section{Direct products, free products, free versions}

In this section we study the behavior of growth invariants under various product operations. The idea is that we have linearisation formulae of type
$$\Phi(A\times B)=\Phi(A)+\Phi(B)$$
where $\Phi$ is some suitably chosen series, equivalent to $B,S$, and $\times$ are various product operations. As for $+$, this will be the usual sum or product of series.

\subsection{Direct products}

Given two pairs $(A_i,u_i)$ as in Definition 1.3, we can form their tensor product
$$(A,u)=(A_1\otimes A_2,u_1+u_2)$$
where the Hopf algebra structure of $A$ is given by the fact that the canonical embeddings $A_i\subset A$ are Hopf algebra morphisms, and $u$ is the direct sum of $u_1,u_2$:
$$u=\begin{pmatrix}u_1&0\cr 0&u_2\end{pmatrix}.$$

\begin{theorem}
For $A=A_1\otimes A_2$ we have $S(z)=S_1(z)S_2(z)$.
\end{theorem}

\begin{proof}
We know that $\Gamma$ consists of products $w=v_1\otimes v_2$ with $v_i\in\Gamma_i$ and
$$l(w)=l(v_1)+l(v_2)$$
$$\dim(w)=\dim(v_1)\,\dim(v_2)$$
and we have the following computation:
\begin{eqnarray*}
S(z)
&=&\sum_{v\in\Gamma} \dim(v)^2 z^{l(v)}\cr
&=&\sum_{v_1\in\Gamma_1}\sum_{v_2\in\Gamma_2} \dim(v_1)^2 \dim(v_2)^2 z^{l(v_1) + l(v_2)}\cr
&=&\sum_{v_1\in\Gamma_1} \dim(v_1)^2 z^{l(v_1)}\sum_{v_2\in\Gamma_2} \dim(v_2)^2 z^{l(v_2)}.
\end{eqnarray*}
The last term is $S_1(z)S_2(z)$, and this gives the assertion.
\end{proof}

\subsection{Free products}

Given two pairs $(A_i,u_i)$ as in Definition 1.3, we can form their free product
$$(A,u)=(A_1*A_2,u_1+u_2)$$
where the Hopf algebra structure of $A$ is given by the fact that the canonical embeddings $A_i\subset A$ are Hopf algebra morphisms, and $u$ is the direct sum of $u_1,u_2$:
$$u=\begin{pmatrix}u_1&0\cr 0&u_2\end{pmatrix}.$$

The invariant which linearises $*$, vaguely inspired from Voiculescu's $R$-transform \cite{vdn}, is constructed as follows.

\begin{definition}
The $P$ invariant of a pair $(A,u)$ is given by
$$P(z)=1-\frac{1}{S(z)}$$
where $S$ is the generating series of numbers $s_k$.
\end{definition}

This invariant is designed to linearise $*$, as we will see in next theorem. We should mention here that an even better invariant, that we won't use here, is
$$Q(z)=\left(1+\frac{1}{z}\right)\left(1-\frac{1}{S(z)}\right)$$
with the normalisation factor being there for having simple formulae for simple algebras. For instance for $A={\rm C}^*(F_n)$ we have $Q(z)=2n$, as one can see by direct computation for $n=1$, then by using the result below for the general case.

\begin{theorem}
For $A=A_1*A_2$ we have $P(z)= P_1(z)+P_2(z)$.
\end{theorem}

\begin{proof}
It is known from \cite{w1} that $\Gamma$ consists of $1$ and of products of type
$$w=v_1\otimes \ldots \otimes v_p$$
with elements $v_i\in\Gamma_1-\{1\}$ alternating with elements $v_i\in\Gamma_2-\{1\}$. The length and dimension of $w$ are given by the following formulae.
$$l(w)=\sum_{i=1}^p l(v_i)$$
$$\dim(w)=\prod_{i=1}^p \dim(v_i)$$

We write $S=1+C_1+C_2$, where $C_i$ is the part of the sum coming from products $v$ satisfying $v_1\in\Gamma_i$. Proceeding as in proof of Theorem~3.1, we get the following formulae relating the series $C_i$ and $S_i$.
$$C_1 = (S_1 - 1) (C_2 + 1)$$
$$C_2 = (S_2 - 1) (C_1 + 1)$$
The solutions of these equations are given by
$$C_1 = \frac{S_1S_2-S_2}{S_1+S_2-S_1S_2}$$
$$C_2 = \frac{S_1S_2-S_1}{S_1+S_2-S_1S_2}$$
so we can compute $S$ by making the sum and adding $1$.
$$S =\frac{2S_1S_2-S_1-S_2}{S_1+S_2-S_1S_2}+1 =\frac{S_1S_2}{S_1+S_2-S_1S_2}.$$
Now by taking inverses we get the formula of $P$ in the statement.
\end{proof}

\subsection{Free versions}

The notion of free version is introduced in \cite{b2}, in a somewhat conceptually heavy context. However, the definition is very simple: the free version of $A$ is the Hopf algebra $\tilde{A}\subset{\rm C}^*({\mathbb Z})*A$ having as fundamental corepresentation the matrix
$$zu=\begin{pmatrix}zu_{11}&&zu_{1n}\cr &\ddots&\cr zu_{n1}&&zu_{nn}\end{pmatrix}$$
where $z$ is the canonical unitary in ${\rm C}^*({\mathbb Z})$, and $u$ is the fundamental corepresentation of $A$. In other words, $\tilde{A}$ is the ${\rm C}^*$-algebra generated by entries of $zu$.

The interest in this notion comes from a diagrammatic formulation of Tannakian duality: the free version appears as universal among pairs $(B,v)$ satisfying
$$P(A,u)=P(B,v)$$
where $P$ is the associated planar algebra. For the purposes of this paper, we need free versions for computing the growth of $A_u(n)$: it is pointed out in \cite{b2} that this is the free version of $A_o(n)$, as one can see by comparing universal properties.

The corepresentation $zu$ is not self-adjoint, so we make the following modification.

\begin{definition}
The free version of a pair $(A,u)$ satisfying $u=\bar{u}$ is the pair $(A^+,u^+)$ constructed in the following way: the fundamental corepresentation is the matrix
$$u^+=\begin{pmatrix}zu&0\cr 0&uz^*\end{pmatrix}$$
and $A^+\subset{\rm C}^*({\mathbb Z})*A$ is the ${\rm C}^*$-algebra generated by coefficients of this matrix.
\end{definition}

Observe that $A^+=\tilde{A}$, and that the corepresentation $v=u^+$ satisfies $v=\bar{v}$.

We assume that all free versions we consider are non-degenerate. Here degenerate means $A^+= {\rm C}^*({\mathbb Z})*A$, and this happens precisely when the trivial corepresentation $1$ is contained in some odd tensor power of $u$. See \cite{b2}.

\begin{theorem}
For a non-degenerate free version $(A^+,u^+)$ we have $P^+(z)=2P(z)$.
\end{theorem}

\begin{proof}
The space $\Gamma^+$ consists of $1$ and of products of type
$$w_1=v_1z^{\pm 1}\otimes v_2z^{\pm 1}\otimes\ldots \otimes v_pz^e$$
$$w_2=zv_1\otimes z^{\pm 1}v_2\otimes\ldots \otimes z^{\pm 1}v_pz^e$$
with $v_i\in\Gamma-\{1\}$, with the choice of each $z^{\pm 1}$ depending on parity of lengths of corepresentations $v_i$ appearing at left of it, and with $e\in\{-1,0\}$ being uniquely determined by parity reasons. The dimension and length of each $w=w_j$ being given by
$$l(w)=\sum_{i=1}^p l(v_i), \quad \dim(w)=\prod_{i=1}^p \dim(v_i)$$
we see that $\Gamma^+$, endowed with its metric and the dimension map $\dim:\Gamma^+\to{\mathbb N}$, is isomorphic to the corresponding space for $A*A$. An isomorphism is indeed given by
$$w_1 \to v_1^{(1)}\otimes v_2^{(2)}\otimes v_3^{(1)}\otimes\ldots \otimes v_p^{(f)}$$
$$w_2\to v_1^{(2)}\otimes v_2^{(1)}\otimes v_3^{(2)}\otimes\ldots \otimes v_p^{(f)}$$
where $v^{(i)}$ is obtained from $v$ by embedding $A$ as $i$-th factor of $A*A$, and where $f\in\{ 1,2\}$ depends on the parity of $p$. The assertion follows now from Theorem 3.2.
\end{proof}

\section{Free quantum groups}

We estimate in this section the growth of various free quantum groups. These yield in particular examples of exponential growth, in the sense that $(b_k)$ dominates a sequence of the form $(r^k)$, with $r>1$.

If $(\log b_k)/k$ has a limit $\log(r)$ as $k \to \infty$, we call $r$ the ratio of exponential growth.  Note that it is not a quasi-isometry invariant, unlike exponential growth.

The finitely generated quantum groups to be considered correspond to universal Hopf algebras endowed with the matrix of their canonical generators and its conjugate:
\begin{eqnarray*}
A_o(n)&=&{\rm C}^*\left( u_{ij}\mid
    u = n\times n\mbox{ orthogonal}\right)\cr
A_u(n)&=&{\rm C}^*\left(
    u_{ij}\mid u =n\times n\mbox{ biunitary}\right)\cr
A_s(n)&=&{\rm C}^*\left( u_{ij}\mid
    u=n\times n \mbox{ magic unitary}\right)
\end{eqnarray*}
These algebras were introduced by Wang in \cite{w1}, \cite{w2}. The precise relations between generators are as follows:
\begin{enumerate} 
\item $u$ orthogonal means $u=\bar{u}$ and $u^t=u^{-1}$.

\item $u$ biunitary means $u^*=u^{-1}$ and $u^t=\bar{u}^{-1}$.

\item $u$ magic unitary means that all $u_{ij}$ are projections, and on each row and each column of $u$ these projections are mutually orthogonal, and sum up to $1$.
\end{enumerate}

The Hopf algebra structure is given by the general formulae in Definition 1.2.

It is known that $A_o(2)$ is a twist of $C(SU(2))$, and that for $n\leq 3$ we have $A_s(n)={\rm C}(S_n)$. We won't consider these algebras: for $A_o(2)$ the invariants are the same as those for $C(SU(2))$, known from section 2, and for $C(S_n)$ the whole problematics is trivial.

In the result below, all algebras are endowed with their canonical
corepresentation $u$. We use the following notations: 
\begin{enumerate} 
\item $q_n$ is the biggest root of $q^2-(n-2)q+1=0$. 

\item $r_n$ is the biggest root of $r^3-(2n^2-1)r^2+2(n^2-1)r-2=0$.
\end{enumerate}

We should mention that our methods can be applied as well to more complicated examples of universal Hopf algebras: for instance to the free analogues of the algebras ${\rm C}(S_a\wr S_b)$, discussed in \cite{b3}. In all cases we have exponential growth, and the ratio can be explicitely computed. We plan to come back to these questions in some future work.

\begin{theorem}
We have the following results, regarding the growth invariants of free discrete quantum groups:
\begin{enumerate}
\item $A_o(n)$ with $n\geq 3$ has exponential growth, with ratio $q_{n+2}^2$.

\item $A_u(n)$ with $n\geq 2$ has exponential growth, with ratio $r_n$.

\item $A_s(n)$ with $n\geq 5$ has exponential growth, with ratio $q_n^2$.
\end{enumerate}
\end{theorem}

\begin{proof}
We use several basic facts on free quantum groups, from \cite{b1}, \cite{b2}, \cite{b3}.

We will see that for the above quantum groups, the series $B,S$ sum up to rational functions. In particular, the study of the poles of these functions will allow us to detect the exponential growth, and to compute the corresponding ratio.

(1) The classification of irreducible corepresentations of $A$, explained in \cite{b1}, is quite similar to that of irreducible representations of $SU(2)$. We have an isomorphism $\Gamma\simeq{\mathbb N}$ so that, denoting by $u_k\in\Gamma$ the corepresentation corresponding to $k\in{\mathbb N}$, $u_0$ is the trivial corepresentation, $\dim u_1 = n$ and the fusion rules read $u_ku_1=u_{k-1}+u_{k+1}$. In particular the sphere of radius $k$ equals $\{u_k\}$ and we have $s_k = (\dim u_k)^2$.

As explained in \cite{b1}, the relation $u_ku_1=u_{k-1}+u_{k+1}$ gives $\dim(u_k)=(a^{k+1}-b^{k+1})/(a-b)$, where $a>b$ are the solutions of $X^2-nX+1=0$. Observe that we have $a+b=n$ and $ab=1$, so we have as well $a^2+b^2=n^2-2$ and $(a-b)^2=n^2-4$. This gives: 
\begin{eqnarray*}
S(z)
&=&\sum_{k=0}^\infty\left(\frac{a^{k+1}-b^{k+1}}{a-b}\right)^2z^k\\
&=&\frac{1}{n^2-4}\sum_{k=0}^\infty(a^{2k+2}+b^{2k+2}-2)z^k\\
&=&\frac{1}{n^2-4}\left(\frac{a^2}{1-a^2z}+\frac{b^2}{1-b^2z}-\frac{2}{1-z}\right)\\
&=&\frac{1}{n^2-4}\left(\frac{(n^2-2)-2z}{1-(n^2-2)z+z^2}-\frac{2}{1-z}\right)\\
&=&\frac{1+z}{(1-(n^2-2)z+z^2)(1-z)}
\end{eqnarray*}

When $n>2$, the biggest root of $z^2 - (n^2-2)z + 1$ is greater than $1$, and
  it coincides in fact with $q_{n+2}^2$:
  \begin{eqnarray*}
    q_{n+2}^4 - (n^2-2)q_{n+2}^2+1 &=& (nq_{n+2}-1)^2 - (n^2-2)q_{n+2}^2 + 1 \\
    &=& 2 q_{n+2}^2 - 2n q_{n+2} + 2 = 0.
  \end{eqnarray*}
  Note that it is also easy to compute $\dim(u_k)$ directly from the fusion
  rules. But we need the expression of $S$ for the next computation.

  (2) We first compute $S$ via the $P$ invariant and Theorem 3.3. Denoting by
  $P_o$, $S_o$ the $P$, $S$ invariants of $A_o(n)$, and by $P_u$, $S_u$ those of
  $A_u(n)$, we have:
  \begin{eqnarray*}
    A_u(n)=A_o(n)^+ &\Rightarrow&P_u=2P_o ~~\Rightarrow~~
    \frac{1}{S_u}=\frac{2}{S_o}-1 \\
    &\Rightarrow& S_u = \frac{S_o}{2-S_o} =
    \frac{1+z}{1-(2n^2-1)z+2(n^2-1)z^2-2z^3}.
  \end{eqnarray*}
  We observe then that the biggest real root $r_n$ of
  $r^3-(2n^2-1)r^2+2(n^2-1)r-2$ is greater than $1$, since
  \begin{eqnarray*}
    &&1^3-(2n^2-1)1^2+2(n^2-1)1-2 = -2 < 0 \quad\text{and}\\
    &&\lim_{r\to\infty} r^3-(2n^2-1)r^2+2(n^2-1)r-2 = +\infty.
  \end{eqnarray*}
  Moreover, there is no complex root of modulus greater than $r_n$, since $b_k
  \in \mathbb{R}$ for all $k$.

  (3) The fusion rules for $A_s(n)$ are the same as those for ${\rm C}(SO(3))$.
  For each $k$ there is a unique length $k$ corepresentation, whose dimension
  can be computed recursively by using the fusion rules:
  \begin{displaymath}
    \dim(u_k)=\frac{(q_n+1)q_n^k-(q_n^{-1}+1)q_n^k}{q_n-q_n^{-1}}.
  \end{displaymath}
  We compute the $S$ series from the
  expression of $\dim(u_k)$ above or from the fusion rules as in (1):
  \begin{displaymath}
    S = \frac{(1+z)^2 + 2(n-2)z}{(1-z)((1+z)^2-(n-2)^2z)}.
  \end{displaymath}

This gives the estimate in the statement.
\end{proof}

\section{Random walks}

Let $(A,u)$ be a finitely generated Hopf algebra as in previous sections, and
consider its Haar functional:
\begin{displaymath}
  \int : A\to {\mathbb C}.
\end{displaymath}
This is the unique positive unital linear form on $A$ satisfying a certain left
and right invariance condition. See Woronowicz \cite{wo}.

\begin{definition}
  For a finitely generated Hopf algebra $(A,u)$, the number
  \begin{displaymath}
    p_k= n^{-2k}\int\left(u_{11}+\ldots+u_{nn}\right)^{2k}
  \end{displaymath}
  with $n={\dim}(u)$ is called probability of returning at $1$ after $2k$ steps.
\end{definition}

The basic example is with $A={\rm C}^*(\Gamma)$, where $\Gamma=<g_1,\ldots
,g_n>$ is a finitely generated discrete group.  The integral here is the linear
extension of $g\to \delta(g,1)$, where $\delta$ is the Kronecker symbol. By
using the identifications in section 1, we see that $p_k$ is indeed the
probability of returning at $1$ after $2k$ steps for the symmetric random walk on $\Gamma$:
\begin{eqnarray*}
  p_k &=&n^{-2k}\int (g_1+\ldots +g_n)^{2k}\cr
  &=&n^{-2k}\sum_{i_1\ldots i_{2k}}\delta(g_{i_1}\ldots
  g_{i_{2k}},1)\cr &=&n^{-2k}\#\{(i_1,\ldots ,i_{2k})\mid
  g_{i_1}\ldots g_{i_{2k}}=1\}.
\end{eqnarray*}

In general, computation of $p_k$ is a representation theory problem, the
quantity to be integrated being the character of the corepresentation
$u^{\otimes 2k}$:
\begin{eqnarray*}
  p_k =n^{-2k}\int\chi_u^{2k} = n^{-2k}\int\chi_{u^{\otimes 2k}}.
\end{eqnarray*}

We can use this method for free quantum groups, and using results from
\cite{b1}, \cite{b2}, \cite{b3} and from \cite{vdn} one gets: 
\begin{enumerate}
\item For $A_o(n)$ with $n\geq 3$ we have $\log p_k\simeq -2k\log(n/2)$.
\item For $A_u(n)$ with $n\geq 2$ we have $\log p_k\simeq -2k\log (n/\sqrt{2})$. 
\item For $A_s(n)$ with $n\geq 5$ we have $\log p_k\simeq -2k\log(n/4)$.
\end{enumerate}

Let us study now the asymptotic behaviour of $(p_k)$ for duals of Lie groups,
which presents an interesting connection with the growth estimate of
Theorem~2.1. 

Recall that connected simply connected compact real Lie groups $G$ are endowed
with the direct sum $u$ of their fundamental representations.

\begin{theorem}
  Let $G$ be a connected simply connected compact real Lie group. Then
  $p_k\simeq k^{-d/2}$, where $d$ is the real dimension of $G$.
\end{theorem}

\begin{proof}
  We first apply Weyl's integration formula \cite[\S 6, cor.~2]{bo9}:
  \begin{eqnarray*}
    p_k &=& n^{-2k}\int_G\chi_u(g)^{2k} \mathrm{d} g
    = w(G)^{-1} n^{-2k} \int_T\chi_u(t)^{2k} \delta_G(t) \mathrm{d}t
  \end{eqnarray*}
  where $T$ is a maximal torus of $G$, $w(G)$ is the cardinality of the Weyl
  group of $G$, and $\delta_G : T \to \mathbb{C}$ is related to the root system
  $R \subset \hat T$ by
  \begin{eqnarray*}
    \delta_G(t) = \prod_{\alpha\in R} (\langle\alpha,t\rangle-1).
  \end{eqnarray*}
  Put $X = \hat T$ and let $R_u \subset X$ be the family of weights of $u$, with
  multiplicity: it generates the lattice $X$ since the coefficients of $u$
  generate ${\rm C}(G)$.  We have by definition $\chi_u(t) = \sum_{\alpha\in
    R_u} \langle\alpha,t\rangle$ for all $t\in T$. We can write, in terms of the
  Fourier transform $\hat\delta_G : X \to \mathbb{C}$ of $\delta_G$:
  \begin{eqnarray*}
    p_k &=& \frac 1{w(G) n^{2k}} \int_T \sum_{(\alpha_i)\in R_u^{2k}}
    \langle\alpha_1 + \cdots +
    \alpha_{2k}, t\rangle~ \delta_G(t) \mathrm{d}t \\
    &=& \frac 1{w(G) n^{2k}} \sum_{(\alpha_i)\in R_u^{2k}}
    \hat\delta_G(-(\alpha_1 + \cdots +\alpha_{2k})).
  \end{eqnarray*}
  Since there are exactly $n$ elements in $R_u$, this shows that $p_k$ can be
  interpreted as the expectation of $\hat\delta_G$ after $2k$ steps of the
  random walk on $X$ generated by $R_u$: if $(S_k)$ denotes the corresponding
  random process, we have
  \begin{eqnarray*}
    p_k &=& w(G)^{-1} \sum_{\alpha\in X} \hat\delta_G(-\alpha) P(S_{2k} = \alpha).
  \end{eqnarray*}

  Now we use the well-known fact that random walks converge to a Gaussian, after
  a suitable renormalization. More precisely, let us denote by $\tilde S_k =
  k^{-1/2} S_k$ the renormalized process with values in $V =
  X\otimes_{\mathbb{Z}} \mathbb{R}$. Then the sequence of characteristic
  functions
  \begin{displaymath}
    \Phi_{\tilde S_k} : V^* \to \mathbb{C},~
    \xi \mapsto E\left(e^{2\mathrm{i}\pi\langle\xi, \tilde S_k\rangle}\right)
  \end{displaymath}
  converges pointwise to $\xi \mapsto e^{-q_1(\xi)/2}$, where $q_1 : V^* \to
  \mathbb{R}$ is the covariance quadratic form of each step of the random walk
  \cite[thm. 29.5]{bi}. We will only need to know that $q_1$ is positive
  definite, but it is not hard to obtain the following formula.
  \begin{displaymath}
    q_1(\xi) = \frac 1n \sum_{\alpha\in R_u} 4\pi^2\langle\xi, \alpha\rangle^2
  \end{displaymath}
  To deal with the renormalization we will need to identify $V^*$ with the Lie
  algebra of $T$ by putting $\langle\alpha, \exp(\xi)\rangle = e^{2\mathrm{i}\pi
    \langle\alpha,\xi\rangle}$ for all $\alpha\in X$, $\xi\in V^*$, and $T$
  itself with the cross-section
  \begin{displaymath}
    \bar T = \{\xi\in V^* ~|~ \forall\beta\in B~
    |\langle\xi,\beta\rangle| \leq 1/2\}
  \end{displaymath}
  associated with a basis $B$ of the lattice $X$.

  We have then, Fourier-transforming back to $\bar T\subset V^*$:
  \begin{eqnarray}\nonumber
    p_k &=& \frac 1{w(G)} \sum_{\alpha\in X} \int_{\bar T} \delta_G(\exp \xi)~
    e^{2\mathrm{i}\pi\langle\alpha,\xi\rangle}~ P(S_{2k} = \alpha)~ \mathrm{d}\xi \\ \nonumber
    &=&  \frac 1{w(G)(2k)^{r/2}} \sum_{\alpha\in X/\sqrt{2k}} \int_{\sqrt{2k} \bar T}
    \delta_G(\exp (\xi/\sqrt{2k}))~
    e^{2\mathrm{i}\pi\langle\alpha,\xi\rangle}~ P(S_{2k} = \sqrt{2k}\alpha)~ \mathrm{d}\xi \\
    \label{eq_asympt_walk}
    &=& \frac 1{w(G)(2k)^{r/2}} \int_{V^*} 1\!\!1_{\sqrt{2k}\bar T}~
    \delta_G(\exp(\xi/\sqrt{2k}))~
    \Phi_{\tilde S_{2k}}(\xi)~ \mathrm{d}\xi.
  \end{eqnarray}
  Let $R_+\subset R$ be a choice of positive roots and put $s = \#R$.  Using the
  following classical computation of $\delta_G$:
  \begin{displaymath}
    \delta_G(\exp \xi) = \prod_{\alpha\in R_+} (e^{2\mathrm{i}\pi \langle\alpha,\xi\rangle} - 1)
    (e^{-2\mathrm{i}\pi \langle\alpha,\xi\rangle} - 1) =
    \prod_{\alpha\in R_+} 4\sin^2(\pi\langle\alpha, \xi\rangle)
  \end{displaymath}
  we see that $(2k)^{s/2}\delta_G(\exp(\xi/\sqrt{2k}))$ converges to the
  following product of quadratic forms.
  \begin{displaymath}
    q_2(\xi) = \Big(\prod_{\alpha\in R_+} 2\pi\langle \alpha,\xi\rangle\Big)^2
  \end{displaymath}
  By the lemma below the bounded convergence theorem applies to
  (\ref{eq_asympt_walk}) and gives, according to the relation $d = r + s$:
  \begin{displaymath}
    \lim_{k\to\infty} (2k)^{d/2}~ p_k = w(G)^{-1} \int_{V^*}
    q_2(\xi) e^{-q_1(\xi)/2} \mathrm{d}\xi \in \left]0,+\infty\right[.
  \end{displaymath}
  This yields the estimate in the statement.
\end{proof}

\begin{lemma}
  The bounded convergence theorem applies to the sequence of
  integrals~(\ref{eq_asympt_walk}) in the proof of the above theorem.
\end{lemma}

\begin{proof}
  By using the estimate $|\sin(t)| \leq |t|$ and notations from previous proof
  we get
  \begin{displaymath}
    0 \leq k^{s/2} \delta_G(\exp(\xi/\sqrt{k}))
    \leq k^{s/2} q_2(\xi/\sqrt{k}) = q_2(\xi).
  \end{displaymath}
  On the other hand we have, by denoting by $S = S_1$ the step of the random
  walk:
  \begin{displaymath}
    \Phi_{\tilde S_{k}} (\xi) = \Phi_S (\xi/\sqrt{k})^{k}.
  \end{displaymath}
  Hence it suffices to prove that $|\Phi_S(\zeta)| \leq e^{-q(\zeta)}$ for all
  $\zeta\in\bar T$ and for some positive definite quadratic form $q : V^* \to
  \mathbb{R}$. This will indeed yield for all $k$:
  \begin{displaymath}
    0\leq |k^{s/2}~ 1\!\!1_{\sqrt{k}\bar T}~ \delta_G(\exp(\xi/\sqrt{k}))~
    \Phi_{\tilde S_{k}}(\xi)| \leq q_2(\xi) e^{- q(\xi)}
  \end{displaymath}
  and the right-hand side is integrable over $\xi\in V^*$.

  By definition of the random walk under consideration, we have
  \begin{eqnarray*}
    &&\Phi_S(\zeta) = \frac 1n \sum_{\alpha\in R_u}
    e^{2\mathrm{i}\pi\langle\zeta,\alpha\rangle} \text{~~~hence} \\
    &&|\Phi_S(\zeta)|^2 = \frac 1{n^2} \sum_{(\alpha,\alpha')\in R_u^2}
    \cos 2\pi\langle\zeta,\alpha'-\alpha\rangle.
  \end{eqnarray*}

  Since we chose $u$ to be the direct sum of the irreducible representations
  associated to the fundamental weights $(\varpi_i)$, the set $R_u$ contains all
  elements $w \varpi$ with $w$ in the Weyl group of $(G,T)$ and $\varpi$ a fundamental
  weight. Now it is an exercise about root systems to check that the set of
  differences $w\varpi-w'\varpi'$ of such elements contains a basis of $X$, which we
  can assume to be $B$. Then for all $\beta\in B$ there exists $(\alpha,\alpha')
  \in R_u^2$ such that $\alpha'-\alpha=\beta$.  Choosing one such pair in
  $R_u^2$ for each $\beta$, and using the evident upper bound $1$ for the other
  terms, we get
  \begin{eqnarray*}
    &&|\Phi_S(\zeta)|^2 \leq 1 + \frac 1{n^2} \sum_{\beta\in B}
    (\cos (2\pi\langle\zeta,\beta\rangle) - 1).
  \end{eqnarray*}
  For every $\beta\in B$ and $\zeta\in\bar T$ we have
  $|2\pi\langle\zeta,\beta\rangle| \leq \pi$. Let us use the inequalities
  $\cos(t) \leq 1-t^2/\pi^2$ for $|t|\leq\pi$ and $1 + t \leq e^t$:
  \begin{displaymath}
    |\Phi_S(\zeta)|^2 \leq 1 - \frac 4{n^2} \sum_{\beta\in B}
    \langle\zeta,\beta\rangle^2
    \leq \exp\Big(-\frac 4{n^2} \sum_{\beta\in B} \langle\zeta,\beta\rangle^2\Big).
  \end{displaymath}
  Hence we can take $q(\zeta) = \frac 2{n^2} \sum \langle\zeta,\beta\rangle^2$,
  which is positive definite since $B$ is a basis.
\end{proof}

Now getting back to discrete groups, a number of estimates relating growth and
random walks are known to hold, see the book \cite{vsc}. In terms of numbers
$b_k$, $p_k$ constructed as above, these estimates make sense as conjectures,
for any discrete quantum group.

We would like in particular to point out the following statement, which in the
case of discrete groups is a key ingredient for Gromov's fundamental result in
\cite{gr}:

\begin{conjecture}
  For a finitely generated quantum group $(A,u)$ we have $b_k\approx k^d$ if and
  only if $p_k\approx k^{-d/2}$.
\end{conjecture}

This statement is verified for discrete groups, by Gromov's result. A second general verification follows from Theorems 2.1 and 6.1 for group duals. Note that these two classes of quantum groups are by no means dual to each other, they rather represent the two ``classical extremes'' of the theory of discrete quantum groups. This is the reason why the proofs of the conjecture for both classes are of very different nature.

As explained in the introduction, a positive answer to the question would be of great interest, in connection with several fundamental problems. We believe also that this conjecture can serve as a good motivation for the study of the world of discrete quantum groups with polynomial growth, which is still largely unknown.

\end{document}